

\def\date{le\ {\the\day}\ \ifcase\month\or janvier\or
{f\'evrier}\or mars\or avril \or mai\or juin\or juillet\or
{ao\^ut}\or septembre\or octobre\or novembre\or {d\'ecembre}\fi
\ {\oldstyle\the\year}}

\let\sskip=\smallskip
\let\bskip=\bigskip

\let\cline=\centerline

\def\a{\alpha} \def\b{\beta}
\def\g{\gamma}

\def\ve{\varepsilon}

\def\t{\theta}

\def\lda{\lambda}
\def\r{\rho}

\def\s{\sigma}

\def\vp{\varphi}

\def\D{\Delta}

\def\Si{\Sigma}

\font\tenbb=msym10

\font\sevenbb=msym7
\font\fivebb=msym5

\newfam\bbfam
\textfont\bbfam=\tenbb \scriptfont\bbfam=\sevenbb
\scriptscriptfont\bbfam=\fivebb
\def\bb{\fam\bbfam}

\def\W{{\bb W}}

\def\C{{\bb C}}

\def\Ff{{\bb F}}

\def\un{{\rm 1\mkern-4mu  l }}

\font\titre=cmbx12

\font\got=eufm10

\def\part{\partial}

\def \wtilde {\widetilde}
\def \oli {\overline}

\def\ra{\rightarrow}
\def\longra{\longrightarrow}
\def\la{\leftarrow}

\def\sbs{\subset}

\def\W{{\cal W}}

\def\and{\mathop{\rm and}\nolimits}

\def\Diff{\mathop{\rm Diff}\nolimits}

\def\Hom{\mathop{\rm Hom}\nolimits}

\def\id{\mathop{\rm id}\nolimits}
\def\Id{\mathop{\rm Id}\nolimits}

\def\Ker{\mathop{\rm Ker}\nolimits}

\def\Spec{\mathop{\rm Spec}\nolimits}

\catcode`\@=11
\def\Eqalign#1{\null\,\vcenter{\openup\jot\m@th\ialign{
\strut\hfil$\displaystyle{##}$&$\displaystyle{{}##}$\hfil
&&\quad\strut\hfil$\displaystyle{##}$&$\displaystyle{{}##}$
\hfil\crcr#1\crcr}}\,} \catcode`\@=12

\catcode`\@=11
\def\displaylinesno #1{\displ@y\halign{
\hbox to\displaywidth{$\@lign\hfil\displaystyle##\hfil$}&
\llap{$##$}\crcr#1\crcr}}

\def\ldisplaylinesno #1{\displ@y\halign{
\hbox to\displaywidth{$\@lign\hfil\displaystyle##\hfil$}&
\kern-\displaywidth\rlap{$##$}
\tabskip\displaywidth\crcr#1\crcr}}
\catcode`\@=12

\def\buildrel#1\over#2{\mathrel{
\mathop{\kern 0pt#2}\limits^{#1}}}

\def\build#1_#2^#3{\mathrel{
\mathop{\kern 0pt#1}\limits_{#2}^{#3}}}

\def\hfl#1#2{\smash{\mathop{\hbox to 6mm{\rightarrowfill}}
\limits^{\scriptstyle#1}_{\scriptstyle#2}}}

\def\hfr#1#2{\smash{\mathop{\hbox to 6mm{\leftarrowfill}}
\limits^{\scriptstyle#1}_{\scriptstyle#2}}}

\def\vfl#1#2{\llap{$\scriptstyle #1$}\left\downarrow
\vbox to 3mm{}\right.\rlap{$\scriptstyle #2$}}

\def\vfll#1#2{\llap{$\scriptstyle #1$}\left\uparrow
\vbox to 3mm{}\right.\rlap{$\scriptstyle #2$}}

\def\diagram#1{\def\normalbaselines{\baselineskip=0pt
\lineskip=10pt\lineskiplimit=1pt}   \matrix{#1}}

\def\up#1{\raise 1ex\hbox{\sevenrm#1}}

\def\signed#1 (#2){{\unskip\nobreak\hfil\penalty 50
\hskip 2em\null\nobreak\hfil\sl#1\/ \rm(#2)
\parfillskip=0pt\finalhyphendemerits=0\par}}

\def\TeX{T\kern-.1667em\lower.5ex\hbox{E}\kern-.125em X}

\def\lsim{ {\raise -3mm \hbox{$<$} \atop \raise 2mm
\hbox{$\sim$}} }

\def\gsim{ {\raise -3mm \hbox{$>$} \atop \raise 2mm
\hbox{$\sim$}} }

\def\frac#1#2{\mathop{\scriptstyle#1\over\scriptstyle#2}\nolimits}

\def\fnote#1{\advance\noteno by 1\footnote{$^{\the\noteno}$}
{\eightpoint #1}}

\def\boxit#1#2{\setbox1=\hbox{\kern#1{#2}\kern#1}%
\dimen1=\ht1 \advance\dimen1 by #1 \dimen2=\dp1 \advance\dimen2 by
#1
\setbox1=\hbox{\vrule height\dimen1 depth\dimen2\box1\vrule}%
\setbox1=\vbox{\hrule\box1\hrule}%
\advance\dimen1 by .4pt \ht1=\dimen1
\advance\dimen2 by .4pt \dp1=\dimen2 \box1\relax}

\def\cube{
\raise 1 mm \hbox { $\boxit{3pt}{}$}
}

\def\cqfd{\unskip\kern 6pt\penalty 500
\raise -2pt\hbox{\vrule\vbox to10pt{\hrule width 4pt
\vfill\hrule}\vrule}\par}

\def\dstar {\displaystyle ({\raise- 2mm \hbox
{$*$} \atop \raise 2mm \hbox {$*$}})}

\def\ref #1#2{
\smallskip\parindent=1,0cm
\item{\hbox to\parindent{\enskip\lbrack{#1}\rbrack\hfill}}{#2} }

\def\choose#1#2{\mathop{\scriptstyle#1\choose\scriptstyle#2}\nolimits}

\def\adots{\mathinner{\mkern2mu\raise1pt\hbox{.}
\mkern3mu\raise4pt\hbox{.}\mkern1mu\raise7pt\hbox{.}}}

\def\pegal{\mathrel{\vbox{\hsize=9pt\hrule\kern1pt
\centerline {$\circ$}\kern.6pt\hrule}}}
\magnification=1200
\cline {\titre SACKLER LECTURES}

\bskip
\cline{\bf  Joseph BERNSTEIN}
\bskip

\cline {\titre Lecture 1. Tensor categories.}
\bskip
{\bf 1.1. Monoidal categories}.
\medskip
(1.1.1)   A {\it monoidal category} is a collection of the following data:

(i) A category $C$;

(ii) A functor $F : C \times C \ra C$ (notation: $ (X,Y) \mapsto X
\otimes Y$);

(iii) An associativity constraint,  i.e. a system of  isomorphisms
$$
a_{X,Y,Z} : (X \otimes Y) \otimes
Z \ \widetilde \ra \  X \otimes (Y \otimes Z)\; \hbox{for any}\;
X, Y, Z \in C.
$$
The system of isomorphisms  $a_{X,Y,Z}$ should be functorial in
$X,Y,Z$ and satisfy the pentagon identity (PI):
$$
\leqno{(PI)}\matrix {
\qquad & (X \otimes(Y \otimes Z))\otimes T &  & \longra & & X \otimes ((Y
\otimes Z)\otimes T) \cr \cr
& \downarrow  &&&&\downarrow \cr
\cr
& ((X \otimes Y)\otimes Z) \otimes T & \ra & (X \otimes Y) \otimes
(Z\otimes T) &
\la & X\otimes (Y \otimes (Z \otimes T)) \cr}.
$$
The notion of a monoidal category is a generalization of the notion of a
monoid.
\medskip

(1.1.2) We usually assume that a monoidal category has a unit object. A
{\it unit
object }  is an object $\un$ together with functorial isomorphisms
$\a_X : X \ \wtilde \ra \ X \otimes \un$ and $\b_X : X \ra \un \otimes X$,
compatible with the associativity constraint.
\bskip
(1.1.3) In the classical theory of monoids the most important example is
that of a group,
i.e.,  a monoid in which all elements are invertible.
 Given an object $X$ in a monoidal category, one can define a notion of
 an inverse object $X^{-1}$, but this notion is rarely useful.
There exists, however,  a weaker notion of the same type which
{\it is} useful.
\medskip
{\bf Definition.} We say that an object $X \in C$ is {\it left rigid} if there
exist an object $ Y \in C$ and morphisms
$$
e : Y \otimes X \ra \un \; {\rm and }\;  i : \un \ra X \otimes Y
$$
such that the following compositions are the identity morphisms ${\rm
id}_X$ and
${\rm id}_Y$, respectively:
$$
\matrix{
X  \ra \un \otimes X\ra X \otimes Y \otimes X  \ra X \otimes \un \ra X \cr
Y  \ra Y \otimes \un\ra Y \otimes X \otimes Y  \ra \un \otimes Y \ra Y }.
$$
\medskip
The object $Y$ and morphisms $e$ and $i$ are defined uniquely up to a canonical
isomorphism. We will denote this object $Y$ by $X^*$ and call it the {\it right
dual} to $X$.
\medskip
We similarly  introduce the notion of a {\it right rigid object} and denote by
${}^*X$  its {\it left dual} (so $({}^*X)^* = X$).
\medskip
The monoidal category $C$ is called a {\it rigid} one if all its objects are
left and right rigid.
\medskip
(1.1.4) {\bf Example}. Let $G$ be a group. Then the category ${\rm Rep}(G)$ of
finite dimensional representations of $G$ is rigid. In this case $X^* = {}^*X$
is the dual representation.
\medskip
(1.1.5) It is easy to see that in a rigid monoidal category we have a
canonical isomorphism $(X \otimes Y)^* = Y^* \otimes X^*$. This
implies that the functor $X \mapsto X^{**}$ has a canonical
structure of a monoidal functor, i.e.,  it is equipped with
  a canonical functorial
 isomorphism $(X \otimes Y)^{**} \approx X^{**} \otimes Y^{**}$.
\medskip
(1.1.6) Another important classical notion connected with monoids is
the notion of a ring. In case of monoidal categories this corresponds
 to the additive monoidal category $C$, such that the multiplication
functor $F$ is biadditive. An example of such category: the category of
representations of a group
( see example (1.1.4)), or of representations of a quantum group, which we
will discuss in Lecture 2.

Note that we may have an abelian monoidal category in which all
objects are rigid. This has no analogues in the classical case.
\medskip
(1.1.7) Let $C$ be an additive monoidal category. We can define a
generalization of the notion of a module
over a ring: a module category over $C$. Namely, such a category consists
of an additive category $
\cal M$, a  biadditive multiplication
functor $H:C \times \cal M \ra \cal M$ and an associativity constraint
$b_{X,Y,M}:(X \otimes Y) \otimes M \widetilde \ra X \otimes
(Y \otimes M) $, satisfying the pentagon identity.
\bskip
{\bf 1.2. Tensor categories}. \bskip
(1.2.1) Let $C$ be a monoidal category. A {\it symmetry constraint } $ S $
is a collection of isomorphisms
$S = \{S_{XY} : X \otimes Y \, \widetilde \ra \, Y \otimes X$ for any $X,
Y\in C \}$
satisfying two hexagon axioms (H1) and (H2):
$$
\leqno{(H1)}
\matrix {
\qquad & X \otimes(Y \otimes Z) & \ra & (X \otimes Y) \otimes Z &
\hfl {S^+} {} & Z \otimes (X \otimes Y) \cr
\cr
& \downarrow
\Id \otimes
S^+&&&&\downarrow \cr
\cr
& X \otimes (Z \otimes Y) & \ra & (X \otimes Z) \otimes Y &
\hfl {S^+
\otimes \Id}
{} & (Z \otimes X) \otimes Y \cr
}$$
where $S^+ = S $.
\sskip
  Axiom $(H2)$ is obtained from $(H1)$ by replacing family $S^+$ with
the family $S^-$ defined by $ S^-_{XY} = (S_{YX})^{-1}  $.
\medskip
(1.2.2) Until about 10 years ago it was always
assumed that $ S^- = S^+$, i.e.,  $S_{YX} \cdot S_{XY} = {\rm id} $
(nowadays such categories are called {\it symmetric} ones).
For example, the category ${\rm Rep}
(G)$ of representations of a group $G$ is symmetric.

On the other hand, the category of representations of a quantum group
satisfies a symmetry constraint but is not symmetric.
\medskip
{\bf Definition.} A {\it tensor category} is an abelian
rigid monoidal category $C$ equipped with a symmetry
constraint $S$.
\medskip
(1.2.3) A symmetry constraint $S$ determines for any pair of objects $X,Y \in
C$ a functorial automorphism  $S_{YX} \circ S_{XY}$ of the object $X \otimes
Y$. If the category $C$ is not symmetric, this automorphism $S_{YX} \circ
S_{XY}$  is nontrivial. In many interesting examples there is an additional
structure  which rigidifies this isomorphism.
\medskip
{\bf Definition.} A {\it balancing} on a tensor category $(C,S)$ is an
automorphism $t$ of the identity functor on $C$ such that
$$
S_{YX} \circ S_{XY}   = t_{X \otimes Y} \circ t^{-1}_X \circ t^{-1}_Y
$$
and
$$
(t_{X^*})  = (t_X)^* .
$$
A tensor category $C$ equipped with a balancing $t$ is called a {\it balanced
tensor category}.
\medskip
(1.2.4) For a balanced tensor category we can identify the two dual objects,
$X^*$ and ${}^*X$. Namely, for any tensor category $C$ we have a canonical
morphism $\a_X : X \ra X^{**}$, given by the composition
$$
\a_X : X \ra X \otimes \un \ra X \otimes X^* \otimes X^{**} \ra X^* \otimes
X \otimes X^{**} \ra \un \otimes X^{**} \ra X^{**}\ .
$$
This morphism is not a morphism of monoidal functors. However, if we consider
the morphism $\b_X = \a_X \circ \,  t_X^{-1} : X \ra X^{**}$,  then $\b_X$ is
an
isomorphism of monoidal functors $\Id $ and $ X\mapsto X^{**}$ from $ C $
to  $ C$.
\medskip  Using this isomorphism
we will identify $X$ with $X^{**}$; therefore, ${}^*X$ with $X^*$.
\bskip
{\bf 1.3. Invariants of knots}.\medskip
(1.3.1) Let $C$ be a tensor category, $X$ and $Y$ two
objects in $C$.
 By permuting $X$ with $Y$ we turn $X \otimes Y$ into
$Y \otimes X$. These two objects are isomorphic, but we have to choose
between two natural isomorphisms
$ S^+ = S_{XY}$ and $S^- = S_{YX}^{-1}$.

Informally speaking, we may say that an isomorphism between the product
$X \otimes Y$ and the permuted product $Y \otimes X$ depends on how
$X$ and $Y$ moved pass each other: if $X$ leaped \lq\lq over" $Y$
we will use $S^+$, if $X$ crawled
\lq\lq under" $Y$ we will use $S^-$.
\medskip
If we have $n$ objects $X_1, \ldots, X_n$, then for any permutation
$\s$ of indices $\{ 1, \ldots , n \}$ we can consider the object
$Z_\s = X_{\s (1)} \otimes \ldots \otimes X_{\s (n)}$.
All these objects are isomorphic, but the choice of an isomorphism
depends on the way the objects $X_i$ leap over or crawl under each other.
\medskip
(1.3.2) There is a geometric notion which generalizes that of
permutation and takes into account the over/under relation between
permuted objects. This is the notion of {\it Artin's braids}.

Let us recall that a {\it braid } $b$ acting on a set $ I = \{
1,\ldots,n\}$ and realizing the permutation $\s$ is a continuous family
of   imbeddings $b_u: I \to  \bb C$, which starts with the
identity imbedding at $u = 0$ and ends with the imbedding defined
by $\s$ at $u = 1$. The set $B_n$ of isotopy classes of such
braids has a natural group structure.
It is called the {\it braid group} of order $n$.
\medskip
{\bf   Claim.} Let $b \in B_n$ be a braid  acting on the set $\{
1,\ldots,n\}$ and realizing the permutation $\s$. Then for objects $X_1,
\ldots   , X_n \in C$ the braid $b$ induces a well-defined
isomorphism  $\g_b : Z_e \,  \widetilde \ra \, Z_\s$.
Moreover, $\g_{b_1b_2} = \g_{b_1} \circ \g_{b_2}$.

Thus, starting from a purely algebraic object --- a tensor category
--- we can construct a representation of such a geometric object
as the braid group $B_n$.
\bskip
(1.3.3) Suppose we have fixed a balanced tensor category $(C,S,t)$.
Then in a way similar to (1.3.2) we can construct an algebraic
representation of another geometric object.

Namely, let $L$ be a {\it  link}, i.e.,  a collection of knotted
oriented circles in ${\bb R}^3$. Suppose that to each circle $\a$ we
have assigned a {\it colouring}, which is an object $X_\a \in C$.
We want to define a weight $w =  w(L,\{X_\a \})$.

Let us choose a generic oriented plane $ M $ in ${\bb R}^3$ and a
generic  linear function $y$ on $M$. We will consider the projection of our
link $L$ on the plane M and study its intersection with horizontal lines.

For any horizontal straight line $\lda$  (given by  the equation $ y = \lda $)
which intersects the projection of $L$ in the general position we consider
an object $X_\lda \in C$, given by $X_\lda = W_{\nu_1}
\otimes W_{\nu_2} \otimes \ldots \otimes W_{\nu_k}$. Here $\nu_1, \ldots ,
\nu_k$ are intersection points of the line $\lda$ with the projection of $L$
(the order of these points is determined by
the orientation of $\lda$ defined by
the  equation $ y = \lda $ and the orientation of $M$).
For every point $\nu_i$  the factor  $W_{\nu_i}$ equals
either $X_\a$ or $X_\a^*$, where $\a$ is the circle which passes through
$\nu_i$, and we choose $X_\a$ if the circle goes up at this point and
$X_\a^*$ if it goes down.
\medskip
(1.3.4) Let us see what will happen with the object $X_\lda$ when we move
the level $\lda$ from $-\infty$ to $\infty$.
It is clear that locally we will only encounter movements of the
following  four types:

I. Two neighboring objects $X = W_{\nu_i}$ and $Y = W_{\nu _{i+1}}$
interchange, so that $X$ leaps over $ Y$.

II. Two neighboring objects $X$ and $Y$ interchange so that $ X$
crawls under $Y$.

III. Two neighboring objects $X$ and $X^*$ collide and dissappear.

IV. Two neighboring objects $X^*$ and $X$ are born out of
thin air.
\medskip
\noindent In other words, we can choose a sequence
$\lda_1 < \lda_2 < \ldots < \lda_N$ such that
\medskip
1) $\lda_1 \ll 0$, so $X_{\lda_1} = \un$;
\smallskip
2) $\lda_N \gg 0$, so $X_{\lda_N} = \un$;
\sskip
3) Passing from $\lda_i$ to $\lda_{i+1}$ we only encounter one of the above
four
simple moves.
\medskip
In each of the above cases I -- IV define  morphisms
$m_i:X_{\lda_i} \ra X_{\lda_{i+1}}$ by setting, respectively:
$$
\matrix{{\rm I)}&S^+  : X \otimes Y \ra Y \otimes X;\qquad\qquad&{\rm
III)}&e:X \otimes X^* \ra \un ;\cr
{\rm II)}& S^-: X \otimes Y \ra Y \otimes X;\qquad\qquad&{\rm IV)}&i : \un
\ra X^* \otimes X.}
$$
\medskip
Now, for any  $i < j$ we define a morphism $m_{ij} : X_{\lda_i} \ra X_{\lda_j}$
as the composition
of morphisms $m_i, \ldots , m_{j-1} $;
$ m_{ij} = m_{j-1}\circ\ldots\circ m_i$.

In particular, we have defined a morphism $m_{1n} : \un \ra \un$.

If we assume that ${\rm End} (\un) = \bb C$, then this morphism
$m_{1n}$ is a number. This number, which we denote by $w(L,\{X_\a \};M,y)$
is the {\it weight} we wanted to describe.
\medskip
(1.3.5) According to definition (1.3.4) the weight $w$ depends on the choice
of  an
oriented  plane $M$ and
the choice of the horizontal direction --- an ordinate $y$ --- on $M$.
In fact, the
axioms of a balanced tensor category imply  that, to a large extent,
this weight only depends on $L$ and the colourings $X_\a$.
\sskip
   In order to get invariants independent of $M$ and $y$ we have to pass
to a {\it framed} link. So let us consider a framed link $L$, i.e.,  a link $L$
together with a field $f$ of nonzero normal vectors on $L$.
   First assume that our framed link $L$ is in a {\it good position}, so that
we can choose an oriented plane $M$ in such a way, that
the frame $f$ on $L$ is everywhere positively transversal to $M$.
Then one can show that the weight $w(L,\{ X_\a \};M,y)$ does not depend on
$M$ and $y$.

Now it is easy to see that any framed link $(L,f)$ is isotopic to a framed link
$(L',f')$ in a good position. Axioms of a rigid balanced tensor category
imply then that the resulting weight $w(L',f')$ only depends on the isotopy
class
of the coloured framed link $(L,f,\{ X_\a \})$.
\sskip
Thus we see, that a balanced tensor category gives a way to produce
invariants of framed links.
\medskip
(1.3.6) {\bf Remark}. The construction of sec. 1.3.5 produces, in fact, a
 representation of an algebraic object -- the category of {\it tangles}.
\medskip
(1.3.7) How to produce invariants of links independent of
colouring? One way to do it is  to assign the same colouring to all circles.
A physical interpretation of our picture
suggests, however,  a better way to do it. Namely, it suggests
to consider the sum of the weights over all possible colourings.

We will see in Lecture 3 that this procedure
works quite well when the category $C$ has finite number of
simple objects $X_\a$. Summing over all possible colourings of the
circles of our link L \lq\lq with" these simple objects produces an invariant
of
framed links $w(L,f)$, which only depends on our balanced category $C$.
\bskip

\cline {\titre Lecture 2. Quantum Groups.}
\bskip
Quantum groups were introduced independently by Drinfeld and Jimbo around 1985.
The term \lq\lq quantum group"  was coined by Drinfeld who also proposed
their theory. Quantum groups soon became very popular and were further studied
by many people. I will mostly follow Lusztig, Majid and Manin.
 \sskip
I will discuss quantum groups from the point of view of tensor categories.
\medskip
{\bf 2.1.  Definition of quantum groups}.
\medskip
(2.1.1) Let $G$ be a finite group, $A =
\C[G]$ the algebra of functions on $G$ with respect to
the
pointwise multiplication. This algebra has the following
structures :
 \sskip
(i) $A$ is an associative algebra with the unit element
$1$. Thus we have two morphisms
$$
m:A \otimes A \ra A\;  {\rm (multiplication)}\quad \quad  {\rm and}\quad \quad
\eta : \C \ra A\;  {\rm (unit)}.
$$
\sskip
 (ii) The multiplication map $G \times G \ra G$ defines a
morphism $\Delta :A \ra A \otimes A $
(comultiplication)  by the formula $\Delta f(x,y) =
f(xy)$. The imbedding of the identity  $\{e\}\  \hookrightarrow
\ G$ defines a morphism $\ve : A \ra \C$ (counit)   by the formula
$\ve(f) = f(e)$.
\medskip

(2.1.2) {\bf Definition.}
A {\it bialgebra} over $\C$ is a linear space $A$ equipped
with the operations
$$ \matrix{
 m &: & A \otimes A & \ra &A \hfill & {\rm (multiplication)}\hfill &
\eta   & : & \C \hfill& \ra &A \hfill& {\rm (unit)}\hfill  \cr
\Delta &: & A \hfill  & \ra &A \otimes A & {\rm (comultiplication)}&
\ve   & : & A \hfill & \ra &\C  \hfill & {\rm (counit)}\hfill   \cr
}$$
satisfying the following axioms H1, H2 and H3 that connects H1 with H2:
$$
\matrix{\matrix{
{\rm H}1.1:\quad& m\;  \hbox{is associative}; \hfill &&
{\rm H}2.1:\quad& \Delta\; \hbox{is coassociative};\hfill  \cr
{\rm H} 1.2:\quad& \eta\; \hbox{is a unit with respect to}\ m; \hfill&&
{\rm H}2.2:\quad& \ve\; \hbox{is a counit with respect to}\ \Delta ;\cr
&&&&\cr}\cr
\matrix{{\rm H}3.1:\qquad\qquad& \Delta \; \hbox {\rm is a morphism of
algebras};\hfill  \cr
{\rm H}3.2:\qquad\qquad& \ve \; \hbox{is a  morphism of algebras}; \hfill \cr
{\rm H}3.3:\qquad\qquad& \eta \; \hbox{is a morphism of coalgebras}.\hfill }}
$$
These three sets of axioms can be expressed as the commutativity of some
diagrams,
constructed in terms of $m, \Delta ,\eta , \ve $. For example, the axiom
H 3.1 is usually written in a symmetric form
$$
\diagram {
A \otimes A \otimes A \otimes A & \hfr{\Delta \otimes \Delta}{}
&A \otimes A&
 \hfl{m}{} & A \cr
\vfl{}{S_{2,3}} & &&& \vfl{}{\Id} \cr
A \otimes A \otimes A \otimes A & \hfl {}{m \otimes m} &
A \otimes A &\hfr{\Delta}{} &A. \cr
}
$$
It is equivalent to the  statement that $m$ is a morphism of coalgebras.
\medskip

(2.1.3) {\bf Definition.}  A bialgebra $A$ is called a {\it Hopf algebra} if
there exists an {\it antipode morphism}
$inv : A\to  A$ such that the two diagrams represented by the following figure
for $i = inv \otimes \id$ and
for $i = \id \otimes inv$  are commutative:
$$
\diagram {
A \otimes A & &\hfl {i}{}&  & A \otimes A       &                  \cr
\vfll{}{\D} & &           &    & \vfl {}{m}      &                  \cr
A           & \hfl{\ve}{} & \C & \hfl {\eta}{}  & A &                  \cr}
$$

It is easy to check that the antipode morphism, if exists, is uniquely
defined. It reverses $m$ and $\D$.
\medskip
(2.1.4) {\bf Example.} Let $A$ be a finite dimensional commutative
semisimple algebra over $\C$. Then it can be realized as $A = \C[G]$, where
$G$ is the finite set $G = \Spec A$. A comultiplication morphism $\D : A \ra A
\otimes A$ defines a multiplication map $G \times G \ra G$. If $A$ is a
bialgebra, this map defines on $G$  the structure of a monoid with the unit
given by $\ve$.  If $A$ is a Hopf algebra, then $G$ is a group.
\medskip

(2.1.5) More generally, let $A$ be a commutative finitely generated algebra
over $\C$. For simplicity assume that it does not have nilpotent elements.
Then it can be realized as the algebra of regular functions on an algebraic
variety $G$. To define on $A$ the  structure of a bialgebra is the same as to
equip $G$ with the structure of an algebraic monoid.
The condition that $A$ is a
Hopf algebra means that $G$ is an  affine algebraic group.
\medskip

(2.1.6) Let $(A,m,\D,\ve,\eta)$ be a finite dimensional Hopf algebra.

Consider the dual vector space $A^*$ and the adjoint operators
$ m^*, \D^* ,\eta^* , \ve^*$.
        Then $ (A^*, \D^*,m^*,\eta^*,\ve^*)$ is also a Hopf algebra.
It is called the {\it dual Hopf algebra} to $A$.
\sskip
For  infinite  dimensional Hopf algebras this definition does not work since
 the comultiplication morphism $m^*$ is not well-defined.  The situation,
however,  can often be mended by considering an appropriate completion
of $A^* \otimes A^*$, or passing to a subalgebra $U \sbs A^*$ on which
$m^*$ is defined.
\medskip

(2.1.7) In all examples of Hopf algebras $A$ we have considered, $A$
is either commutative (e.g. $A = \C [G])$, or
cocommutative (e.g. $A = \C[G]^*$).
A natural question is whether there exist {\it natural} examples of Hopf
algebras, which are neither commutative nor cocommutative.
\medskip

What Drinfeld and Jimbo have discovered is a family of such Hopf algebras
(Drinfeld called them quantum groups). I will formulate the result of Drinfeld
and Jimbo in a form closer to Manin's description.
\medskip

{\bf Statement}. {\it Let $G$ be a simple algebraic group over
$\C,\ A = \C[G]$ the Hopf algebra of regular functions on $G$. Then there
exists a nontrivial family
of Hopf algebras (quantum groups) $A_q$ parametrized
by $q \in \C^*$ such that} $A_1 = A$.

\medskip

(2.1.8) Let us describe in detail the case of $G = SL(2,\C)$.
In this case the algebra $A$ is generated by four indeterminates
$a,\ b,\ c,\ d$ which commute and satisfy the relation $ad - bc = 1$.

In order to describe the
comultiplication $\D : A \ra A \otimes A$ we introduce a matrix
$$
Y = \pmatrix
{  a & b \cr c & d \cr } \in {\rm Mat} (2,A)\ .
$$
Using the
natural imbeddings $i',i'' : A \ra A \otimes A,\ \left ( i'(x) = x
\otimes 1,\ i''(x) = 1 \otimes x \right )$, we can define morphism
 $\D : A \ra A \otimes A$ by the formula
$$
\D(Y) = i'(Y) \cdot i''(Y)\ ,\leqno (*)
$$
which is an equation considered in ${\rm Mat}\;  (2, A \otimes A)$.
Condition $(*)$ is a
shorthand for equations
$$
\D a = a \otimes a  + b \otimes c,\;
\D b = a \otimes b + b \otimes d, \quad \quad \hbox {\rm and so on.}
$$
\sskip
Now we can describe the quantum group $SL_q(2)$. It is given by an algebra
$A_q$, generated by elements $(a,b,c,d),$ satisfying the following relations
$$
\matrix{
ab &= q^{-1} ba &\qquad ac &= q^{-1}ca &\qquad bc &= cb \cr
cd &= q^{-1}dc & bd &= q^{-1}d b &ad - da &= (q^{-1}-q)bc \cr
&&&ad - q^{-1}bc = 1&& \cr}\leqno{(**)}
$$
\sskip

The comultiplication $\D$ is given by the same condition $(*)$ as above.
\medskip

In Manin's Montreal notes there is a beautiful proof which shows that $A_q$
is a Hopf algebra and explains the origin of
conditions $(**)$.
\medskip

{\bf 2.2. Representations of quantum groups}.
\medskip
(2.2.1) Let $\rho$ be a representation of a finite group $G$ in a vector
space $V$. It is possible to define $\rho$ in terms of the coalgebra $A=\C [G]$
as  a morphism $\rho : V \ra A \otimes V$ which satisfies the following axioms
(roughly speaking the axioms of a $G$-module with reversed arrows in the
commutative diagrams that  define it):
\medskip
 R1) $\rho$ is coassociative, i.e.,  the two natural morphisms $V \ra A
\otimes A
\otimes V$ coincide,

R2) the counit $\ve$ acts on $V$ as the
identity.

\noindent This example is a model for the following general notion.
\medskip

(2.2.2) {\bf Definition.} Let $A$ be a coalgebra. A {\it comodule} over
$A$ is a morphism $\rho : V \ra A \otimes V$, satisfying axioms $R1,\ R2$.
\medskip
(2.2.3) If $V$ is a comodule over $A$, then the through map $V \hfl{\rho}{}
A \otimes V
\hfl {\vp}{} V$ for any element $\vp \in A^*$ defines
an operator \sskip
$\rho (\vp) : V \ra V$.
\sskip
   The axioms $R1,\ R2$ are equivalent to the condition that $V$ is an
$A^*$-module.
\medskip

When $A$ is  finite dimensional this functor $ A$-comod $\to$ $A^*$-mod
 gives an equivalence between the categories of $A$-comodules and
$A^*$-modules. If $A$ is infinite dimensional  this functor is fully faithful,
but it is not an equivalence of categories.
\medskip

(2.2.4) {\bf Example.} Let $V = \C^n$. Then a
morphism $\rho : \C^n \ra A \otimes
\C^n$ defines a matrix $Y \in {\rm Mat} (n,A)\; $ by $\rho (e_i) = \Si y_{ij}
\otimes e_j$. The condition that $\rho$ defines a comodule structure on
$\C^n$ is equivalent to the condition that the matrix $Y$ is {\it
multiplicative}, i.e.,   $\D(Y) = i'(Y) \otimes i''(Y)$ (cf. 2.1.8).
\medskip

(2.2.5) It is natural to define a {\it finite quantum group} as a finite
dimensional Hopf algebra $A$. Let us also define an {\it algebraic quantum
group} as a Hopf algebra $A$ such that
\medskip

(i) $A$ is a finitely generated algebra,

(ii) $A$ is a union of finite dimensional $A$-comodules.
\medskip
\noindent By analogy with the classical situation we can formulate two
problems:
\medskip
I)  How to describe finite quantum groups.

II)  How to describe algebraic quantum groups.
\medskip

(2.2.6) Let us discuss Problem II in a particular case of
quantum groups $G_q$ which are deformations of a classical simple group $G$. It
turns out that for groups $G$ other than $SL(2)$ it is difficult to give
an explicit description of the algebra $        A_q$.
The reason for this is that the algebra $A = \C[G]$ is already quite
complicated.
 (Do you know the description of this algebra for the group $G$ of
type $F_4$ ?)

  In order to describe this algebra somehow, we usually pass to the Lie
algebra $\hbox {\got g}=Lie (G)$ and to the enveloping algebra
$ U (\hbox {\got g})$. This algebra $U (\hbox {\got g})$ lies in $A^*$ and
inherits the structure of a Hopf algebra from $A^*$. While $A^*$ is,
clearly, too big, the algebra $U (\hbox {\got g})$ is quite manageable.
\sskip
Similarly, in the quantum case one usually describes instead of the complicated
Hopf algebra $A_q$ a simpler Hopf algebra, $U_q \sbs A^*_q$. Then, in terms of
the algebra $U_q$, it is usually not difficult to reconstruct the Hopf algebra
$A_q$.
\medskip

(2.2.7) {\bf Example.} $G_q  = SL_q(2).$

Let $I \sbs A_q$ be an ideal generated by $b$ and $c$, $S = A_q/I =
\C[a,d]/(ad-1)$ (see (2.1.8)). Then $I$ is a Hopf ideal in $A_q$ and $S$
is a Hopf algebra isomorphic to the Hopf algebra of regular functions on
the algebraic group $H = \C^*$.
   Informally, this means that our quantum group $G_q$ contains $H$ as a
subgroup.
\sskip
In particular, the group $H$ acts on the algebra $A_q$ from the left and from
the right.
\medskip

Consider elements $E,F,K \in A^*_q$ defined as follows. Identify $A/I^2$ with
$S \oplus S\,b \oplus S\, c$ and define $E,F,K : A/I^2 \ra \C$ by
$$
K(P+Qb + Rc) = P(q);\;
E(P+Qb + Rc) = Q(q);\;
F(P+Qb + Rc) = R(1) .
$$
\medskip

{\bf Theorem.} {\it Let $U_q$ be a
subalgebra of $A^*_q$ generated by $E,F,K^{\pm
1}$. Then $U_q$  is a Hopf algebra given by the following relations}
$$
\matrix{KE = q^2EK;\; KF = q^{-2}FK;\;
[E,F] = {K - K^{-1} \over q - q^{-1}}\cr
\D\, K = K \otimes K;\;
\D\,E = E \otimes 1 + K \otimes E;\;
\D\, F= F \otimes K^{-1} + 1 \otimes F \ .}\leqno{(L)}
$$

This is Lusztig's form of generators for $U_q$.
\medskip

(2.2.8) For the general simple group $G $   Lusztig has described the
algebra $ U_q$ as a Hopf algebra, generated by generators $E_i, F_i,
K_i^{\pm1}$ satisfying some relations similar to $(L)$. \medskip
Representations of the quantum group $G_q$ can be realized as some special
$U_q$-modules.
\medskip

{\bf 2.3. Representations of a Hopf algebra as a monoidal
category.}
\medskip

(2.3.1) Let $(\rho,V)$ and $(\s,E)$ be two comodules over a Hopf algebra
$A$. Then
using the multiplication morphism $m :  A \otimes A \ra A$ we can define an
$A$-comodule structure on the space $V \otimes E$ by setting
$$
V \otimes E \hfl{\rho\otimes\sigma}{}\ (V \otimes A) \otimes (E \otimes A)
\hfl{S_{2, 3}}{}\ V \otimes E \otimes A \otimes A \hfl {m}{} V \otimes E
\otimes
A \ .
$$

Let $ {\rm Rep}(A)$ be the category of
finite dimensional $A$-comodules. Then the
tensor product described above defines on ${\rm Rep}(A)$ the
structure of a monoidal
category. Using the antipode one can show that this category is rigid.
\medskip

The forgetful functor $(\rho,V) \mapsto V$ defines a monoidal functor
$H : {\rm Rep}(A)  \ra {\rm Vec}$, where ${\rm Vec}$
is the category of finite dimensional vector spaces.
\medskip

(2.3.2) The following theorem, due to Majid, is philosophically significant,
since it clarifies the relation between Hopf algebras and tensor categories.
\medskip

{\bf Definition.} Let $R$ be a monoidal category. A {\it fiber functor } is a
functor $H : R \ra {\rm Vec}$, together with a functorial isomorphism
$$\eta_{_{X,Y}} : H(X \otimes Y) \ \widetilde \ra \ H(X) \otimes H(Y)\ ,$$
which is compatible with the associativity constraints in the categories $R$
and ${\rm Vec}$.
\medskip

(2.3.3) {\bf Theorem}. {\it Let $R$ be an abelian
rigid monoidal category and $H
: R \ra {\rm Vec}$ a fiber functor. Then there
exists a Hopf algebra $A$ and an
equivalence of monoidal categories $R \simeq {\rm Rep}(A)$
under which $H$ becomes
the forgetful functor. The Hopf algebra $A$ is uniquely defined up to
a canonical isomorphism}.
\medskip

(2.3.4) Let $R$ be an abelian rigid monoidal category. Suppose it has a fiber
functor $H : R \ra {\rm Vec}$. Every such
functor leads to a Hopf algebra $A = A(H)$.
For nonisomorphic functors these Hopf algebras are not isomorphic but it is
clear that one should consider them as equivalent since they describe
essentially the same algebraic object.
\sskip
The corresponding notion of equivalence for Hopf algebras was introduced by
Drinfeld. Namely, consider an element $Q \in (A \otimes A)^*$. Then using the
left and right $A \otimes A$-comodule structures on $A \otimes A$ we define
the operators $\rho_L(Q)$ and $\rho_R(Q)$ on $A \otimes A$. Note that these
operators only depend on the coalgebra structure on $A$.
\medskip

Now, assume that the operator $\rho_R (Q)$ is invertible and consider a new
algebra structure on $A$ given by the multiplication
operator $m_Q : A \otimes A \ra A$, where $m_Q = m \circ \rho_L(Q) \circ
\rho_R(Q)^{-1}:A \otimes A \ra A$. It is easy to see that $m_Q$ is always
a morphism of coalgebras.
\medskip

Observe that the
associativity condition imposes rather complicated restrictions on $Q$.
Drinfeld avoids these complications by introducing the notion of a {\it
quasi-Hopf
algebra}, a generalization of the notion of
Hopf algebra.) Drinfeld
calls the new Hopf algebra $A_Q  = (A,m_Q,\D)$ {\it gauge equivalent} to the
original Hopf algebra $A =(A,m,\D)$.
\medskip

Since, as a coalgebra, $A_Q = A$, it is clear that we have a natural
(identity) equivalence of
categories $I : {\rm Rep}(A) \simeq {\rm Rep}(A_Q)$.
This equivalence can be extended to
an equivalence of monoidal categories. Namely, given two $A$-comodules
$(\rho,V)$ and $(\s,E)$ we define an isomorphism $\a_{\rho,J\sigma}: I(\rho)
\otimes I(\s) \ \widetilde \ra \  I(\r \otimes \s)$ using an operator $(\rho
\otimes \s)(Q) \in {\rm Aut}(V \otimes E)$,
where both $A_Q$-comodules $I(\r) \otimes
I(\s)$ and $I(\r \otimes \s)$ are realized in the same vector space $V \otimes
E$ and the automorphism $(\r \otimes \s)(Q)$ defines an isomorphism between
them.
\medskip

(2.3.5) After the original Drinfeld-Jimbo
construction of a 1-parametric family of deformations of $A=\C[G]$ several
people noticed that for every simple group $G$ there actually exists a
family of deformations of the Hopf algebra $A$, which depends on a large number
of parameters and the  Hopf algebras in the family are nonisomorphic. (For
example for $G = SL(n,\C)$ this family of Hopf algebras depends on $\simeq
n^2/2 $ parameters.)  Later Drinfeld showed that  every
deformation of $A$ only depends on one parameter (say, the original
Drinfeld-Jimbo's one) if considered up to a gauge equivalence.
\medskip

(2.3.6) Given a Hopf algebra $A$, we can consider a new Hopf algebra
$A^\circ$ with the opposite multiplication $m^\circ$ and the same
comultiplication $\D$. When $A$ is a deformation of the algebra $ \C[G] $,
 Drinfeld's results imply that this new Hopf algebra $A^\circ$ is gauge
equivalent to $A$.
\medskip

Let $R \in (A \otimes A )^*$ be an element which defines this equivalence
(Drinfeld calls it the {\it universal $R$-matrix}). It is easy to see that
$R$ defines a symmetry constraint $S_{XY} : X \otimes Y \ \widetilde \ra
 \ Y \otimes X$ on the category ${\rm Rep}(A)$. Thus,
${\rm Rep} (A)$ is a tensor category.
\medskip

Actually, Drinfeld has also constructed a balancing on this category. So
${\rm Rep}(A)$ is a balanced tensor category.
\medskip

(2.3.7) Another natural possibility suggested by Theorem 2.3.3 is
that there might exist an abelian rigid monoidal category $R$ which does not
have any
fiber functor. Then $R$ is not directly related to any Hopf algebra, though
intuitively it represents a mathematical object of the same nature. We will
see natural examples of such categories in 2.4.
\medskip

(2.3.8) The category $R = {\rm Rep}(A)$ is usually studied from the
dual point of view. Namely, one fixes a Hopf subalgebra $U \sbs A^*$ and
considers the category $R$ of finite dimensional $U$-modules with tensor
product given by $\D$, i.e.,
$$
U \otimes (V \otimes E) \hfl{\D}{}\, U \otimes U \otimes (V \otimes E) \ra
(U \otimes V) \otimes (U \otimes E) \ra V \otimes E \ .
$$
Drinfeld, Jimbo and Lusztig only work with this dual picture, and
only study the Hopf algebra $ U_q$ with practically no reference to the
algebra $A_q$.
They also use slightly different subalgebras $U$ inside the algebra $A^*_q$.
\medskip

(2.3.9) Quantum groups $U_q$ studied by Lusztig are deformations of the
algebra $U_1= U(\hbox {\got g})$, which is the  enveloping algebra
of a simple Lie algebra {\got g}. Lusztig showed that any finite
dimensional representation $V $ of $U(\hbox {\got g})$ admits a
deformation $V_q$ which is a representation of the algebra
$U_q$ in the same vector space $V$.
\medskip

For a generic $q$ this construction allows us to describe all
$U_q$-modules. When $q$ is a root of 1 the situation is much more
interesting and complicated. For example, in this case the category
${\rm Rep} (G_q) $, as well as the category of
$U_q$-modules, are not semisimple.
\medskip

{\bf 2.4. Finite quantum groups.}
\medskip

(2.4.1) Let $G$ be a simple algebraic group $A = \C[G]$. Fix a prime
number $l$ and consider a
quantum group $A_q$, where $q = \root l\of{1}$. Then
there exists a natural Frobenius morphism of Hopf algebras $Fr : A \ra
A_q$, which can be interpreted as a homomorphism of quantum groups
$Fr^*: G_q \ra G$.
\medskip
(2.4.2) {\bf Example}. Consider $G = SL(2)$. Then $Fr$ is given by
$$ \pmatrix {
a & b  \cr
c & d  \cr
}
\qquad \longmapsto \qquad
\pmatrix {
a^{l} & b^{l}  \cr
c^{l} & d^{l}  \cr }.
$$
\medskip

(2.4.3) Let $H_q$ be the kernel of the homomorphism of quantum
groups $Fr^* : G_q \ra G$.
This kernel is a finite quantum group, given by the finite
dimensional Hopf algebra $H_q = A_q/J$, where $J$ is the ideal
generated by $Fr \bigl(\Ker (\ve:A \ra \C) \bigr )$.
For example, if $G = SL(2)$ the ideal $J$ is generated by
$\{b^{l}, c^{l}, a^{l}-1,d^{l}-1 \}$.

The quantum group $H_q$ is the quantum analogue of the finite group
$G(\Ff_l)$. This analogy is explored in detail in Lusztig's works.
\medskip
(2.4.4) The category $R_q = {\rm Rep}(H_q)$ is usually not semisimple.
It turns
out that one can modify it (replace by an appropriate quotient category) so
that the resulting category $\oli R_q$ is a semisimple tensor category with
a finite number of simple objects. Let us describe how to
do it in case of $SL(2)$.
\medskip

(2.4.5) The group  $G = SL(2,\C)$  has a series of irreducible representations
 $V_i,\ i = 0, 1, 2,\ldots, $ with $  \dim V_i = i + 1$. For any $q$
we will denote also by $V_i$ the deformed representations of the quantum group
$G_q$.
\medskip

Let $q = \root l\of{1}$. Let $V$ be an $H_q$-comodule. Then $V$ is an
$A_q$-comodule, and using formulas in (2.2.7) and (2.2.3)
 we can define operators $E, K, F : V \ra V$.
\medskip

We say that the $H_q$-comodule $V$ is {\it negligible} if it is free as
$\C[E]/(E^{l})$-module.
\medskip

Let us define the category $\oli R_q$ as the quotient of the category $R_q$
modulo all negligible objects. One can
show that this notion is well-defined and
the resulting category $\oli R_q$ is an abelian semisimple category with
isomorphism classes of simple objects given by $V_0,\ldots,V_{l-2}$ (these
modules are simple objects in $R_q$).
   Since the tensor product of a negligible
object by any other object is negligible, $\oli R_q$ inherits the structure of
a balanced tensor category from the category $R_q$.
\medskip

(2.4.6) Thus, for $q = \root l\of{1}$ we have described a semisimple balanced
tensor category with a finite number of simple objects. One can construct
similar categories starting with any simple group $G$.
\medskip
One can check that these categories do not have fiber functors, so they do not
admit a direct description in terms of Hopf algebras as in 2.3; we have,
however, described them using Hopf algebras $A_q$ and $H_q$.
\medskip
These categories are the sources of nontrivial invariants of knots
described in Lecture 1 and of topological field theories which we
are going to discuss in Lecture 3.
\bskip

\cline {\titre Lecture 3. Topological (quantum) field theories.}
\bskip
{\bf 3.1. Topological field theories ( TFT)}.
\medskip
TFT emerged in physics in the study of conformal field theories. They were
first
explicitly studied by E. Witten. The mathematical framework for TFT is
described in  Atiyah's paper on TFT.
\medskip
(3.1.1) Let us fix
terminology: a {\it manifold} is a compact oriented manifold
$M$ with
boundary. Its boundary $\part M$ will be considered with the induced
orientation. We say that $M$ is closed if $\part M = \emptyset$.
\medskip
We will describe TFT of
dimension $(d+1)$ (a physical interpretation: we have
$d$ space variables and $1$ time variable).
\medskip
We will usually use the following mnemonic rule: we denote the
$(d+1)$-manifolds (i.e manifolds of dimension $d+1$) by $N$, $d$-manifolds
by $M$,
$ (d-1)$-manifolds by
$L$, and so on.
\medskip
(3.1.2) We are going to describe a topological field theory $W$ of
dimension $(d+1)$. This will be done on several levels.
\medskip
{\bf Definition of TFT on
level 0.} A TFT $W$ is a multiplicative invariant
of closed $(d+1)$-manifolds.

In other words, $ W $ assigns to every closed
$(d+1)$-manifold $N$ a number $w(N)$ such that $w(N \cup N') = w(N) \cdot
w(N')$.
\medskip
(3.1.3) {\bf Definition of TFT on level 1.}
\sskip
A TFT $W$ assigns to every closed $d$-manifold $M$ a finite dimensional
vector space $W(M)$ and to every bordism $N$ between two $d$-manifolds $M_1$
and $M_2$ it assigns an operator $w(N) : W(M_1) \ra W(M_2)$.
\medskip
(3.1.4) In definition (3.1.3) we assume that the correspondence $M \mapsto
W(M)$ is functorial with respect to diffeomorphisms of $d$-manifolds.
\medskip
(3.1.5) We also assume that this correspondence is multiplicative.
This means that in
addition to the functor $M \mapsto W(M)$ we have a functorial
isomorphism
$\g_{M,M'} : W(M \cup M') \, \widetilde \ra \, W(M) \otimes W(M')$.
\medskip
(3.1.6) A {\it bordism} between two closed $d$-manifolds $M_1$ and
$M_2$ is a $(d+1)$-manifold $N$ equipped with an isomorphism $\part N \simeq
M^*_1 \cup M_2$ (here $M^*_1$ is the manifold $M_1$ with the opposite
orientation).
\medskip
(3.1.7) The data, described in (3.1.3), i.e.,  the functor $W$, the
system of functorial isomorphisms $\g_{M,M'}$ and the system of operators
$w(N)$, should satisfy the following requirements.
\medskip
$(i)$ $(W,\g)$ is a monoidal functor compatible with the
natural symmetry constraint.
In other words, $W(\emptyset) = \C$ and the system of
isomorphisms $\g$ is compatible with the  natural
associativity constraint, with the unit
object, and with the natural symmetry constraint.
For example, the last condition
just means that the natural diffeomorphism $M \cup M' \simeq  M' \cup M$
corresponds to the standard isomorphism $W(M) \otimes W(M') \simeq   W(M')
\otimes W(M)$.
\medskip
$(ii)$ {\it Functoriality.} Diffeomorphic
bordisms give the same operator.
\medskip
$(iii)$ {\it Composition}. Given three closed $d$-manifolds
$M_1,\ M_2,\ M_3$ and two bordisms: $N'$ between $M_1$ and $M_2$, and $N''$
between $M_2$ and $M_3$, we can glue $N'$ with $N''$ and get a bordism $N$
between
$M_1$ and $M_3$ (notation: $N = N'' \circ N'$ ).
We require that
$$w(N) = w(N'') \circ w(N')\ .$$
\sskip
$(iv)$ {\it Cylinder.} Let $C(M) = [0,1] \times M$ be the cylinder
bordism between $M$ and $M$. Then $w \bigl (C(M) \bigr ) = \Id_{W(M)}$.
\medskip
$(v)$ {\it Multiplicativity.} If $N$ is a bordism between $M_1$ and
$M_2$ and $N'$ a bordism between $M'_1$ and $M'_2$, then
$$w(N \cup N') = w(N) \otimes w(N').
$$
\sskip
(3.1.8) {\bf Comment 1.} The axioms of TFT imply the following homotopy
invariance property.
\medskip
{\bf Lemma.} Let $\vp_t : M \ra M'$ be a smooth family of
diffeomorphisms. Then locally the morphism
$W(\vp_t) : W(M) \ra W(M') $ does not depend on $t$.
\medskip
This shows, that TFT $W$ defines a representation of the mapping
class group  $Cl(M) =
\Diff(M) /  \Diff^0(M)$ in the vector space $W(M)$.
\medskip
(3.1.9.) {\bf Comment 2.} TFT $ W$ defines an invariant of closed
$(d+1)$-manifolds. Indeed, every such manifold $N$ can be considered
as a bordism between $M_1 = M_2 = \emptyset$.
Hence $w(N) \in \Hom (\C,\C) = \C$. This
shows the connection of descriptions on level 1 and level 0.
\medskip
Note that usually TFT can be uniquely reconstructed from this invariant of
$(d+1)$-manifolds.
\medskip
(3.1.10) {\bf Comment 3.} Let $N$ be a $(d+1)$-manifold and $M=\part N$. We
can interpret $N$ as a bordism between the empty manifold $\emptyset$ and $M$.
Then the operator $w(N) \in \Hom \bigl(\C,W(M)\bigr )$ is nothing but a
vector $w(N) \in W(M)$. Similarly, if $\part N=M^*$, we can interpret $w(N)$
as a functional on $W(M)$.
\medskip
(3.1.11) {\bf Comment 4}. Let $N$ be a closed $(d+1)$-manifold. Suppose we
have cut it into two pieces, $N_1$ and $N_2$, with the common boundary $M$.
Then $w(N_1)$ is a vector in $W(M)$ and $w(N_2)$ is a functional on $W(M)$.
The composition property implies that
$$ \langle w(N_2),\ w (N_1) \rangle = w(N)\ .$$
\medskip
This shows, that for a TFT $W$ the corresponding invariant of
$(d+1)$-manifolds has a
local character in the following sense: the number $w(N)$ can
be computed by cutting $N$ in pieces. An important and highly nontrivial
feature of topological field theories is that the answer does not depend on
the nature of the cut.
\medskip
(3.1.12) {\bf Comment 5.} For a closed $d$-manifold $M$ we can interpret the
cylinder $N = C(M)$ as a bordism between $M \cup M^*$ and the empty manifold
$\emptyset$. This defines a canonical pairing $w(N) : W(M) \otimes W(M^*)
\ra \C$. The
axioms of TFT imply that this pairing is nondegenerate; so
it defines a canonical isomorphism $W(M^*) \simeq  W(M)^*$.
\medskip Similarly we have a canonical morphism $w(N) : \C \ra W(M) \otimes
W(M^*)$, which induces the same isomorphism $W(M^*) \simeq W(M)^*$.
\medskip
{\bf 3.2.  Fusion Algebra}.
\medskip
We will be mostly interested in the case $d=2$. But first let us look at the
case $d=1$.
\medskip
(3.2.1) Fix a TFT $W$ of dimension $(1+1)$.
\sskip
Consider an (oriented) unit circle $S$ and set $A = W(S)$. This is a
well-defined vector space, since the group $\Diff^+(S)$ of orientation
preserving diffeomorphisms of $S$ acts trivially on $A$.
\sskip
We have a canonical (up to homotopy) diffeomorphism $\t : S \ra S^*$. It
defines an isomorphism $\t : W(S) \ra W(S^*) \simeq W(S)^*$, i.e.,  a
bilinear form $\t$ on $A$. It is easy to see that $\t$ is symmetric and
nondegenerate.
\medskip
(3.2.2) Let $N$ be a 2-sphere with 3 disjoint discs removed ($N$ is usually
called {\it pants}). We consider $N$ as a bordism between $S \cup S$ and
$S$. Then the operator $w(N)$ defines a
multiplication
$$m : A \otimes A \ra A\ .$$
The axioms of TFT imply that $m$ is associative (this follows from the
geometric fact that bordisms $N \circ (N \circ N)$ and
$(N \circ N) \circ N$ are diffeomorphic).
\medskip
It turns out that $m$ is also commutative. In order to see this
it suffices to  consider a rotation $\vp$ of pants through
180${}^\circ$, so that circles $S_1$ and $S_2$ interchange and
$S_3$ turns by 180${}^\circ$.
\medskip
If we take the disc $D$ with $S=\partial D$, then it is easy to see that the
corresponding element $w(D) \in W(S)=A$ is a unit of the algebra $A$.
\medskip
(3.2.3) We can interpret the element $w(D)$ as a functional $\eta : A \ra
\C$. Then, clearly, $\t(a,b) = \eta(ab)$.
\medskip
(3.2.4) It is easy to check that the algebra $A$ equipped with the form $\t$
(or equivalently with the functional $\eta$) allows us to uniquely reconstruct
TFT $W$.
\medskip
In examples the algebra $A$ is usually semisimple, i.e.,  isomorphic to $\build
\oplus_{i=1}^k \C$.
\medskip
(3.2.5) {\bf Exercise.} $\dim A = w(T^2)$, where $T^2$ is the torus
$S \times S$.
 \medskip
(3.2.6) Now consider a $(2+1)$-dimensional TFT $ W$. Starting with it
we can produce a $(1+1)$-dimensional TFT $ V$ by
$V(L) = W(S \times L)$, where $S$ is the standard circle. As we saw,
such theory $V$ can be described by a commutative algebra
$A = V(S) = W(S \times S)$. This algebra $A$ is called the {\it fusion
algebra} of TFT  $ W$.
\medskip
{\bf 3.3.  $\bf (2+1)$-dimensional theories.}
\medskip
(3.3.1) How to give examples of $(2+1)$- dimensional topological field
theories, and hence construct invariants of 3-manifolds?  Let us try
to analyze the situation.
\sskip
To begin with we are mostly interested in a level $0$ description, i.e.,
to every
$3$-manifold $N$ we would like to assign an invariant $w(N)$.
\medskip
We can try to do it by passing to a level 1 theory. Then instead of very
complicated objects --- $3$-manifolds --- we have to study much more manageable
objects --- $2$-manifolds. However, we have to pay for this simplification: now
to a manifold we assign not a number,
but an algebraic object --- a vector space.
\medskip
It turns out, that we can move further in the same direction. Namely, we can
pass to $1$-manifolds, provided we assign to them even more complicated
structures. This naturally leads us to a level 2 description of TFT.
\medskip
(3.3.2) {\bf Definition of TFT on level 2.}
\medskip
I) A TFT $ W$ of dimension $(d+1)$ assigns
to every closed $(d-1)$-manifold
$L$ an abelian category $\W(L)$ of finite type over $\C$ (see below).
\sskip
II) To every bordism $M$ between manifolds $L_1$ and $L_2 $ TFT $ W$
assigns a functor $W(M) : \W(L_1) \ra \W(L_2)$.
\sskip
III) Suppose we are given three $(d-1)$-manifolds $L_1, L_2,L_3$ and
bordisms $M'$ between $L_1$ and $L_2$ and $M''$ between $L_2$ and $L_3$.
Consider a composite bordism $M$ between $L_1$ and $L_3
$. Then TFT $ W$ should  provide an isomorphism
$\a_{M',M''} : W(M) \,\widetilde \ra\, W(M'') \circ W(M')$.
\sskip
IV) Given two bordisms $M$ and  $M'$  between $L_1$ and $L_2$ and
a bordism $N$ between $M$ and $M'$ the TFT $W$ should
assign to $N$ a morphism of functors $w(N) : W(M) \ra W(M')$.
\medskip
(3.3.3) An abelian category of finite type over $\C$ may be defined as
a category equivalent to the category of finite dimensional
$B$-modules for some finite dimensional $\C$-algebra B.

In fact, in all known cases of TFT the algebra B is semisimple, so the
corresponding category is equivalent to the category of vector bundles on
a finite set.

For categories of finite type one naturally defines their tensor product.

(3.3.4) We assume that the correspondence $L \mapsto \W(L) $ in (3.3.2) I)
is functorial with respect to diffeomorphisms of $L$.
This means that to every diffeomorphism $\vp $
it assigns a functor $W(\vp )$, and to every
composition of diffeomorphisms $\vp \circ \psi$ it assigns an
isomorphism of functors
$W(\vp \circ \psi) \widetilde \ra W(\vp) \circ W(\psi)$.
These isomorphisms should satisfy some version of the pentagon identity.

Similarly, in (3.3.2) II) we assume that the correspondence $M \mapsto W(M)$ is
functorial with respect to diffeomorphisms of $M$.

We also assume that the correspondence $W: M \mapsto W(M)$ is multiplicative.
\medskip
(3.3.5) The data described in 3.3.2 should satisfy many compatibility
conditions. For example, for a composition of three bordisms $M_1, \ M_2,\ M_3$
there should exist an identity of the type of the  pentagon identity.
\medskip
 When we attempt to list all these conditions, we end up with many
pages of definitions and verifications of their compatibilities.
By working with an appropriate notion of a stalk of manifolds
one can reduce this description to a manageable size.
\medskip
 We will not try to list all these conditions and give a rigorous
definition of level 2 TFT. Instead let us fix a $(2+1)$-dimensional
TFT  $W$ and try to describe what structures
lie behind this notion.
\medskip
(3.3.6) Consider the abelian category $C = \W(L)$, where $L = S$ is the
standard circle. It turns out that the requirement that $\W(S)$
functorially depends on $S$ and the homotopy invariance property
imply that the category $C$ has an additional structure, namely
an automorphism of the identity functor
$t : {\rm Id_C} \ra {\rm Id_C}.$
\medskip
Indeed, consider the family of rotations
$\vp_u : S \ra S\,,\ u \in [0,1]$,
where $\vp_u$ is the rotation of $S$ through the angle of $2\pi u$.
Fix an object $X \in C$. Then for the same reasons as in 3.1.8
the family of objects $X_u = W(\vp_u)(X) \in C$ should be locally
constant. Since the equality of objects does not make sense,
this just means that $X_u$ is a local system of objects of the category $C$
 on the segment $[0,1]$. For example, if $C = {\rm Vec}$  this system is the
usual local system on the segment.
This local system defines canonical isomorphisms between all
objects $X_u$. On the other
hand, since the diffeomorphism $\vp_1$ is identity, we have
$X_1 = X$. This
defines a monodromy operator $t : X = X_0 \, \widetilde \ra \, X_1 = X$.

  Later we will use the following  equivalent description of the
 automorphism $t$. Consider the functor $P = W(\vp_{1/2}): C \ra C$
and the isomorphism $p:{\rm Id} \ra P$ described above. Then $P^2$ is equal
(i.e.,  is canonically isomorphic) to the identity functor ${\rm Id}$, so the
morphism $p$ defines an isomorphism $p^2: {\rm Id} \ra {P^2} = {\rm Id}$,
which we
denote by $t$.
\medskip
(3.3.7) Now, consider pants as in (3.2.2) that define the bordism $M$
between $S \cup S$ and $S$. Then the TFT $W$ assigns to $M$ a
functor $F = W(M) : C \times C \ra C$.
This functor defines on $C$ the structure of a monoidal category.

   Namely, we define the associativity constraint as an isomorphism
corresponding to the natural diffeomorphism
$ M \circ (M \circ M) \widetilde \ra (M \circ M) \circ M $.
The pentagon identity follows from the fact that two composite
diffeomorphisms between products of three bordisms are isotopic.
\medskip
(3.3.8) Now, let us define a symmetry constraint $S$ on $C$.
Let $\vp$ be the rotation of pants $M$ through the angle of 180$^{\circ}$
(as in (3.2.2)).
Then
$W(\vp)$ defines an isomorphism
$$
W(\vp) : P(X \otimes Y)\, \widetilde \ra \, P(Y) \otimes P(X),
$$
where $P : C \ra C$ is the functor, corresponding to the rotation of $S$ by
180$^{\circ}$ (see (3.3.6)).
\medskip
Using an isomorphism $p : \Id \ra P$, described in (3.3.6) we can interpret
$W(\vp)$ as an isomorphism
$$
S_{XY} : X \otimes Y \ra Y \otimes X.
$$
The rotation $\vp$ can be extended to a diffeomorphism $\vp : N \circ (N
\circ N) \,
\widetilde \ra \, (N \circ N) \circ N$; this implies that the isomorphism
$W(\vp)$ is compatible with the associativity constraint. Rewriting this in
term
of $S_{XY}$ gives hexagon identities for the symmetry constraint $S$.
\medskip
(3.3.9) In order to carry constructions described in (3.3.6) - (3.3.8) we
only need
part of the data that defines TFT, namely, data I, II, III. So, let us define a
{\it modular functor} $W$ (of dimension (2+1)) as a correspondence
$$
\matrix{
{\rm I.}&\qquad
\{ 1{\rm -manifold}\  L \}\qquad&
\Longrightarrow & {\rm a\  category}\; W(L). \hfill\cr
{\rm II.}&\qquad
\left \{{\rm A\ bordism}\ M \hfill
\atop  {\rm between}\ L_1\ {\rm and}\ L_2\ \right \} \hfill
& \Longrightarrow &\hbox{ a functor}\ W(M): W(L_1) \ra W(L_2).\hfill\cr
{\rm III.}&
\left \{ {\rm Bordisms}\ M'\ {\rm between}\ L_1\ {\rm and}\ L_2
\hfill\atop
{\rm and}\  M'' \ {\rm between}\  L_2 \ {\rm and}\ L_3 \right \}\hfill
&\Longrightarrow &{\rm an\  isomorphism}\hfill
\atop  a_{M',M''}: W(M) \,\widetilde \ra \, W(M'') \circ W(M').\hfill}
$$
\medskip
Then the constructions described in (3.3.6) -- (3.3.8) show that
the modular functor
$W$ defines a balanced tensor category on $C = \W(S)$.
\medskip
(3.3.10) {\bf Remark}. We do not actually need to know values of the modular
functor $W$ on all bordisms.

Let us say that a bordism $M$ is of {\it genus} $0$ if the closed surface
$M$ obtained from $M$ by gluing discs to all circles on its boundary is
diffeomorphic to a sphere.

Let us define a {\it genus $0$ modular functor} as a correspondence I, II, III
only defined when all bordisms are of genus 0. Then as in (3.3.6) -- (3.3.8) we
check that a genus 0 modular functor $W$ defines a balanced tensor category.

It turns out that any balanced tensor category uniquely corresponds to a genus
0
modular functor. In other words, {\it an algebraic notion of balanced
tensor category
is the same as a geometric notion of a genus 0 modular functor}.

This fact,  implicitly contained in Drinfeld's paper on quasi-Hopf
algebras, was explicitly formulated by Deligne in his letter to
Drinfeld (Deligne
used slightly different but equivalent language of punctured algebraic curves).
\end